\documentclass[a4paper]{amsart}
\usepackage{amssymb}
\title[Computer Search for Curves]{Computer search for curves with many points among abelian covers of genus 2 curves}
\author{Karl R\"okaeus}
\address{Karl R{\"o}kaeus \\ Korteweg de Vries Instituut voor Wiskunde \\ Universiteit van
Amsterdam \\ P.O. Box 94248  \\1090 GE Amsterdam \\ The Netherlands}
\email{S.K.F.Rokaeus@uva.nl}
\date{June 23, 2011}

\newcommand{\comment}[1]{}

\DeclareMathOperator{\PP}{\mathbf{P}}
\DeclareMathOperator{\Cl}{Cl}
\DeclareMathOperator{\Gal}{Gal}
\DeclareMathOperator{\NN}{N}
\DeclareMathOperator{\FF}{\mathbf{F}}

\theoremstyle{remark}

\newtheorem*{example}{Example}
\newtheorem*{acknowledgements}{Acknowledgments}

\begin{document}
\maketitle

\emph{Abstract.}\footnote{The author is supported by The Wenner-Gren
Foundations Postdoctoral Grant}
Using class field theory one
associates to each curve $C$ over a finite field, and each subgroup
$G$ of its divisor class group, unramified abelian covers of $C$ whose
genus is determined by the index of $G$. By listing class groups of
curves of small genus one may get examples of curves with many
points; we do this for all curves of genus 2 over the
fields of cardinality 5,7,9,11,13 and 16, giving new
entries for the tables of curves with many points \cite{manyPoints}.

\section*{Background}
 Let $q$ be a prime power and let $g$ be a non-negative integer. We are interested in how many rational points there can be on a
projective, smooth and geometrically irreducible curve of genus $g$, defined
over the finite field $\FF_q$; equivalently, how many rational places there can be in a global function field of genus $g$ with full constant field $\FF_q$. Denote by $\NN_q(g)$
the maximum possible number of rational places in such a field. The current intervals in
which $\NN_q(g)$ in known to lie for $g\leq 50$, together with references, are given in
the tables \cite{manyPoints}. In this report we improve upon the
lower bounds of these intervals, by constructing curves with many points for some genera
over the fields with $q=5,7,9,11,13$ and $16$ elements. We do this by going through unramified abelian covers of all genus 2 curves over these fields, using class field theory. The method of looking at such covers has been used before by various authors, see \cite{vanderGeer}. Still, by doing the systematic search we found some new curves.

The search was done using Magma \cite{magma}, first in terms of curves, using its arithmetic geometry part; then everything was checked also using its function field part. In this report we have formulated everything in terms of function fields.

\subsection*{Notation}
See \cite{Stich} for an introduction to function fields in one variable. By a function field $F/k$ we will always mean a global function field with full constant field $k$, i.e., $k$ is algebraically closed in $F$.
We write $\PP_F$ for the set of places of $F$ and $\Cl(F)$ for its divisor class group of degree zero divisors.

\section*{Description of the method}
\subsection*{Class field theory for global function fields:}
In this section we state the theorems that we use for the construction of the curves; see \cite{NiederXing} for details. Let $F/k$
be a global function field. Let $O$ be
a rational place of $F$ and set $S=\PP_F\setminus\{O\}$. We use $F_O$ to denote the $S$-Hilbert class field of $F$, i.e., the
maximal unramified abelian extension of $F$ (inside a fixed
separable closure) in which $O$ splits completely.
There is an isomorphism
\begin{align*}
\psi\colon\Cl(F) &\to \Gal(F_O/F).
\end{align*}
under which $[P-\deg(P)O]$ maps to $\bigl[\tfrac{F_O/F}{P}\bigr]$,
the Artin symbol of $P$, for any $P\in\PP_F$. Fix a subgroup $G$ of
$\Cl(F)$ of index $d$ and write $F_{O}^G$ for $F_{O}^{\psi(G)}$, the
subfield of $F_O$ fixed by $\psi(G)$. Then $F_O^G$ is a Galois
extension of $F$ with Galois group isomorphic to $\Cl(F)/G$; in particular $[F^G_O:F]=d$. A rational
place $P\in\PP_F$ splits completely in $F^G_O$ if and only if
$\bigl[\tfrac{F_O/F}{P}\bigr]\in\psi(G)$, i.e., if and only
if $[P-O]\in G$.

Therefore, if $F/k$ is a global function field of genus $g$ with
rational places $\{P_1,\dots,P_m\}$ and if $G\subset\Cl(F)$ is a
subgroup of index $d$ then, for every rational place $O$, we get an
unramified extension $F^G_O/F$ of degree $d$, hence of genus
$d(g-1)+1$, with $d\cdot \left|G\cap\{[P_i-O]\}_{i=1}^m\right|$
rational places. Also, $k$ is the full constant field of $F^G_O$.

Note that, for fixed $G$, two different rational places $O$ and $O'$
need not give isomorphic extensions. Although $F^G_O/F$ and $F^G_{O'}/F$
have the same genus and isomorphic Galois groups  they may have a different number of rational places.

We give an example that can be done by hand, a construction of a
genus 6 curve over $\FF_2$ with 10 rational points.
\begin{example}
Let $F$ be the degree 2 extension of $\FF_2(x)$ given by
$$y^2+xy=x^5+x^3+x^2+x.$$
It is of genus 2 and has 4 rational places; one containing $x$ and one containing $1/x$, $P_{0}$ and $P_\infty$; and two containing $x+1$, $P_{1,0}$ and
$P_{1,1}$. Moreover, $F$ has two places of degree 2 so its class number
is $h=10$. This already gives the existence of four covers of degree 10, namely the $S$-Hilbert class fields $F_O$, where $S=\PP_F\setminus\{O\}$ and $O$ is any rational place in $F$. These are all of genus 11, with 10 rational places.

Let $a=[P_{0}-P_\infty]\in\Cl(F)$. Then $a\neq 0$ (because $F$ is not rational), while $2a=0$ because $(x)=2P_{0}-2P_\infty$. Let
$G=\{0,a\}\subset\Cl(F)$. If we choose $O=P_\infty$ then $[P_\infty-O]=0$ and $[P_{0}-O]=a$ lie
in $G$. The extension $F^G_O/F$ therefore has
degree $10/2=5$. It has genus 6 and contains $2\cdot 5=10$ rational
places.

Similarly, if we choose $O=P_{0}$ we also get a genus 6 extension with
10 rational places. However, if we instead choose $O$ to be $P_{1,0}$ or $P_{1,1}$
then $F_O^G$ is a genus 6 extension with only 5 rational places.
\end{example}

\subsection*{Organization of the search}
Let $k$ be a finite field and let $g>1$ be an integer. Let $\mathcal{F}_g$ be a list of global function fields
over $k$ of genus $g$. (In this report we will only consider the case of $\mathcal{F}_2$, the fields of genus 2.) For each function field $F\in\mathcal{F}_g$, do the following:
\begin{itemize}
\item Compute the class group $\Cl(F)$ of $F$.
\item For each rational place $O$ compute the image $I_O$ of the rational places in $\Cl(F)$ under the map $P\mapsto [P-\deg(P)O]\colon \PP_F\to\Cl(F)$.
\item List all subgroups of $\Cl(F)$. For each of these, $G$, compute its index $d_G$ and the cardinality $m_{G,O}$ of its intersection with $I_O$. This gives the existence of an unramified abelian extension of $F$ of degree $d_G$, of genus $d_G(g-1)+1$ and with $d_{G}m_{G,O}$ rational places.
\end{itemize}

\section*{Implementation and results}
Using Magma \cite{magma} we implemented the search with all genus 2 curves as base, over $k=\mathbf{F}_q$ for $q\leq 16$. For $q$ smaller than 9 the listing of genus 2 curves could be done in the straightforward way; for the bigger values of $q$ we used the Magma package \emph{g2twist} \cite{g2twists}. The curves are then represented by Magma as hyperelliptic curves. For each such curve we computed its Jacobian and a list of its subgroups, using the arithmetic geometry part of Magma. All the records we found were then checked using the function field part of Magma, this is also how we describe them in this report.

Below we have recorded all the new curves we found by this method, together with the details required for their construction: First a genus 2 function field $F$. Then a set of rational places $M$ that generate the subgroup of the class group that defines the extension. More precisely, the extension is $F_O^G$ where $O$ is any place in $M$ and $G$ is the subgroup of $\Cl(F)$ generated by $\{[P-O]\}_{P\in M}$. (It can also be constructed as the $S$-Hilbert class field of $F$, with $S=\PP_F\setminus M$.)

The genus 2 fields are always given as a degree 2 extension $F=k(x,y)/k(x)$ defined by some $f\in k[x][y]$. We use the following notation for the rational places of $F$: For $\alpha\in k$ let $Q_\alpha$ be the place of $k(x)$ containing $x-\alpha$. If $Q_\alpha$ ramifies then there is a unique place of $F$ containing $x-\alpha$ which we denote by $P_\alpha$. If $Q_\alpha$ splits there are two such places, one for each $\beta\in k$ such that $f(\alpha,\beta)=0$; it is characterized by containing $y-\beta$, we denote it $P_{\alpha,\beta}$. If there is a unique place of $F$ containing $1/x$ we denote it $P_\infty$.

\subsection*{Curves in characteristic 2}
This is the case where most of the work of constructing curves has been done, in particular using the present method, and for $q=2,4$ and $8$ we didn't find any improvements.

Over $\FF_{16}$ there are $8470$ different genus 2 curves. Among their unramified abelian covers we found 2 new records:
\begin{itemize}
\item \underline{$g=8$, $N=63$.}  $F$: $y^2+(x^2+x)y=\alpha^6x^5+\alpha^{12}x^4+x^3+\alpha^3x^2+\alpha^9x$, where $\alpha$ is a primitive element of $\FF_{16}$ satisfying $\alpha^4+\alpha+1=0$. Places: three ramified, $P_\infty,P_0,P_1$ and the 3 pairs $P_{\alpha^3,\alpha^4},P_{\alpha^3,\alpha^{10}}$, $P_{\alpha^6,\alpha^6},P_{\alpha^6,\alpha^{12}}$ and $P_{\alpha^9,\alpha^3},P_{\alpha^9,\alpha^9}$. Old interval $[62,75]$.

\item \underline{$g=49$, $N=240$.} $F$: $y^2+(x^2+x)y=x^5+x^3+x^2+x$. Places: With $O$ any of $P_\infty,P_0,P_1$, there are two such extensions. One in which these three places and $P_{\alpha^5,\alpha^2}$ and $P_{\alpha^5,\alpha^8}$ splits; and one in which these three places and $P_{\alpha^{10},\alpha}$ and $P_{\alpha^{10},\alpha^4}$ splits. Old interval $[213,286]$.
\end{itemize}

\subsection*{Curves in characteristic 3}
For $q=3$ we didn't find any new curves. For $q=9$ there had been much previous activity, \emph{e.g.,} \cite{vanderGeer}, \cite{vanderGeer2}, \cite{NiederXing5},
\cite{Shabat}. By doing a systematic search we found 3 small improvements of the old table entries. Let $\alpha$
be a primitive element of $\mathbf{F}_9$ satisfying $\alpha^2-\alpha-1=0$.
\begin{itemize}
\item \underline{$g=18$, $N=68$.} $F$: $y^2=x^5 + \alpha^6x^3 + \alpha^6x^2 +\alpha^3 x$. Places: $P_\infty$, $P_0$, $P_1$ and $P_\alpha$.  Old interval $[67,84]$.
\item \underline{$g=32$, $N=93$.} $F$: $y^2=x^5 + \alpha^6x^4 + \alpha^7x^3 + 2x^2 + \alpha^5x + \alpha^2$. Places: $P_\infty$, $P_{0,\alpha}$ and $P_{0,-\alpha}$. Old interval $[92,130]$.
\item \underline{$g=38$, $N=111.$} $F$: $y^2=x^5 + \alpha^7x^3 + 2x^2 + \alpha^2x$. Places: $P_\infty$, $P_0$ and $P_{-\alpha}$. Old interval $[105,149]$.
\end{itemize}

\subsection*{Curves over $\FF_5$}
For $q=5$ there had been some previous work, \emph{e.g.,} \cite{NiederXing2}, \cite{NiederXing3} and \cite{NiederXing4}, in particular using the present method. By doing the systematic search we found 3 improvements of the old table entries:

\begin{itemize}
\item \underline{$g=7$, $N=24$.} $F$: $y^2=x^5+x^2-x$. Places: $P_\infty$, $P_0$, $P_{4,1}$ and $P_{4,4}$. Old interval $[22,26]$.
\item \underline{$g=9$, $N=32$.} $F$: $y^2=x^5-x^3+x$. Places: $P_\infty,P_0$ and either the two places containing $x-1$ or the two places containing $x-4$. Old interval $[26,32]$.
\item \underline{$g=12$, $N=33.$} $F$: $y^2=x^5+x^4+3x^2+1$. Places: $P_\infty$, $P_{0,1},P_{0,4}$ . Old interval $[30,38]$.
\end{itemize}

We give some details about the construction of the genus 9 curve that we found which, together with the Oesterl\'e bound, show that $\NN_5(9)=32$.
\begin{example}
Let $k=\FF_5$ and let $F/k$ be the hyperelliptic field of genus 2
given by
$$y^2=x^5-x^3+x.$$
Take $O$ to be $P_0$, the unique place containing $x$. Magma gives an isomorphism from $\Cl(F)$ to
$\mathbf{Z}/8\oplus\mathbf{Z}/8$ under which the divisor classes
$[P_\infty-O],[P_{4,3}-O]$ and $[P_{4,2}-O]$ map to $(0,4), (0,5)$
and $(0,3)$. They are hence contained in a
subgroup of index 8. Therefore, $F$ has an unramified extension of
degree 8 in which these places and $O$ split completely. This extension
has genus $1+8\cdot(2-1)=9$ and $8\cdot 4=32$ rational places.

We mention that Magma can give explicit equations for the class fields: The degree 8 extension of $F$ in this example has minimal polynomial
\begin{multline*}
T^8+((4x+4)y+3x^4+3x^3+3x^2+3x+3)T^4+\\
+(x^5+2x^4+2x^3+2x^2+2x+1)y+x^8+2x^6+x^5+2x^4+x^3+2x^2+1.
\end{multline*}
\end{example}

\subsection*{Curves over $\mathbf{F}_{7},\FF_{11}$ and $\mathbf{F}_{13}$}
For $q=7,11$ and $13$ there where no lower bounds in the tables for genus greater than 4, so all we had to do to make it to the tables was to produce curves with more than $b(g)/\sqrt{2}$ points, where $b(g)$ is the best known upper bound for $\NN_q(g)$.
(This is the criteria used to decide if a curve  is considered to have many points, see \cite{vanderGeerVlugt}.) The results are given in the table below; we give integers $g$ and $N$ such
that there exists a genus $g$ curve with $N$ points, and
then the interval $[b(g)/\sqrt{2},b(g)]$.

The details of the construction are given in \cite{Rokaeus} and \cite{Rokaeus2}. There, for each entry given in this table, a function field of genus 2 together with all the rest of the information required to prove the existence of the the extension, are given.

\noindent\begin{minipage}{\textwidth} \vspace{.2in}\noindent{\large
New entries for the tables over $\mathbf{F}_7,\FF_{11}$ and $\mathbf{F}_{13}$.}\\[.05in]
\begin{tabular}{|l|l|l||l|l|l||l|l|l|}\hline
$q=7$ &  & &  $q=11$ &  & & $q=13$ & &\\
\hline
$g$ & $N$ & $[\frac{b(g)}{\sqrt{2}},b(g)]$ & $g$ &  $N$ & $[\frac{b(g)}{\sqrt{2}},b(g)]$& $g$ &  $N$ & $[\frac{b(g)}{\sqrt{2}},b(g)]$ \\
\hline
 5 & 24 & [20,28] & 5 & 36 &  [27,38] & 5 & 40 & [32,44]\\
 6 & 25 & [23,32] & 6  & 40 & [32,45] & 6  & 50 & [37,52]\\
 7 & 30 & [26,36] & 7  & 42 & [36,50]& 7  & 48 & [42,58]\\
     8 & 35 & [27,38] & 8 & 42 & [39,55] & 8 & 56 & [45,63]\\
     9 & 32 & [29,41] & 9 & 48 & [42,59]& 9 & 56 & [49,68] \\
     10 & 36 & [32,45] & 10 &54 & [45,63]& 10 & 63 & [51,72]\\
     11 & 40 & [34,48] & 11 &60  & [48,67]& 11 & 70  & [55,77]\\
     12 & 44 & [37,51] & 12 & 66 & [51,72]& 12 & 66 & [58,82]\\
     13 & 48 & [39,54] & 13 & 60 & [55,77]& 13 & 72 & [62,87]  \\
     14 & & & 14 & & & 14 & 65 & [65,91] \\
     15 & & & 15 &70 & [61,85]& 15 & 84 & [68,96]\\
     16 & 45 & [45,63] & 16 & & & 16 & 90 & [72,101] \\
     17 & 64 & [47,66] & 17 &80 & [66,93]& 17 &96 & [75,105]  \\
     18 & 51 & [49,68] & 18 & 85 & [70,98]& 18 & 85 & [78,110] \\
     19 & 54 & [51,71] & 19 &  90 & [73,102]& 19 &  90 & [82,115] \\
     20 & & & 20& 76 & [75,106]& 20& 95 & [85,119]\\
     21 & 60 & [55,77] & 21& 80 & [78,110]& 21& 100 & [88,124]\\
     22 & 63 & [56,79]   &  22   &    &       &22 & 105 & [92,129] \\
     23 &66 & [58,82] & 23 & 88 & [85,119]  & \dots  & &               \\
     24 & & & 24 & 92 & [87,123]  & 29 & 140 & [114,161]\\
     25 & 72& [62,87] & 25 & 96 & [90,127]& 32 & 124 & [124,175]\\
     26  &    &            & 26  &      &            & 33 & 128 & [127,179]\\
    27 & 78 & [66,93]  & 27 & &  & 35 & 136 & [133,187]        \\
         & & &  28 &108 & [98,138] & 36 & 140 & [136,191]  \\
         & & &  \dots &  & & 37 & 144 & [138,195]\\
         & & &  33 & 128 & [111,156] & 41 & 160 & [150,211]\\
         & & &  &  & & 43 & 168 & [155,219] \\
\hline\end{tabular}\end{minipage}
\newline

\begin{acknowledgements}
The author thanks the Wenner-Gren Foundations for
financial support; the Korteweg-de Vries Institute at the University of
Amsterdam for hospitality; and Gerard van der Geer for his
hospitality and for helpful conversations and comments.
\end{acknowledgements}

\comment{\bibliographystyle{amsplain}
\bibliography{biblio}}

\providecommand{\bysame}{\leavevmode\hbox
to3em{\hrulefill}\thinspace}

\end{document}